\documentclass{ifacconf}

\usepackage{graphicx}      
\usepackage{natbib}        

\usepackage{amsmath} 
\usepackage{amssymb}  
\usepackage{epsfig} 
\usepackage{booktabs}
\usepackage{colortbl}
\usepackage{multirow}
\newcommand{\R} {\mathbb{R}}

\newcommand{\La} {\mathcal{L}}
\newcommand{\bbma} {\begin{bmatrix}}
\newcommand{\ebma} {\end{bmatrix}}

\newtheorem{remark}{Remark}

\definecolor{mygray}{gray}{0.9}
\definecolor{mydarkgray}{gray}{0.85}
\definecolor{darkgray}{RGB}{169,169,169}
\definecolor{darkgray176}{RGB}{176,176,176}
\definecolor{gray}{RGB}{128,128,128}
\definecolor{green01270}{RGB}{0,127,0}
\begin{document}
\begin{frontmatter}

\title{Efficient Nonlinear MPC by  Leveraging LPV Embedding and Sequential Quadratic Programming\thanksref{footnoteinfo}} 

\thanks[footnoteinfo]{This work was funded by the German Research Foundation (DFG), project number 419290163.}

\author[First]{Dimitrios S. Karachalios} 
\author[First]{Hossam S. Abbas}

\address[First]{The Institute for Electrical Engineering in Medicine, University of L\"ubeck,  L\"ubeck, Germany (e-mail:dimitrios.karachalios@uni-luebeck.de, h.abbas@uni-luebeck.de).}

\begin{abstract}                
In this paper, we present efficient solutions for the nonlinear program (NLP) associated with nonlinear model predictive control (NMPC) by leveraging the linear parameter-varying (LPV) embedding of nonlinear models and sequential quadratic programming (SQP). The corresponding quadratic program (QP) subproblem is systematically constructed and efficiently updated using the scheduling parameter from the LPV embedding, enabling fast convergence while adapting to the behavior of the controlled system. Furthermore, the approach provides insight into the problem, its connection to SQP, and a clearer understanding of the differences between solving NMPC as an NLP and using the LPV-MPC approach, compared to similar methods in the literature. The efficiency of the proposed approach is demonstrated against state-of-the-art methods, including NLP algorithms, in control benchmarks and practical applications.

\end{abstract}

\begin{keyword}
Linear parametrically varying methodologies, Control of constrained systems, Nonlinear predictive control, Optimal control,  Quadratic and nonlinear programming.
\end{keyword}

\end{frontmatter}

\section{Introduction}

Nowadays, model predictive control (MPC)  is becoming one of the most essential modern control strategies in practice and research \citep{QiBa03, RaMaDi17chapter8}. One strength of MPC is its simple concept, which can accommodate variant classes of control problems with constraints. Its strategy is to compute the control trajectory after solving constrained optimization problems within a reasonable sampling time. Often, nonlinear (NL) MPC (NMPC) is more relevant in practice due to the inevitable emergence of nonlinear models in interesting applications. 
However, implementing NMPC in practice by casting its optimization
as a nonlinear program (NLP) to control high-order nonlinear systems with long prediction horizons and to achieve complex tasks, such as trajectory planning \citep{alcala2020autonomous}, remains a challenge.

A simple and practical approach to mitigating the computational complexity of NMPC is to use the linear parameter-varying (LPV) framework \citep{To10} by introducing the so-called scheduling parameter $p$, which depends on the state of the system and/or the input, and its value is assumed to be measurable (available) at every sampling instant. If the values of the scheduling parameter over the MPC prediction horizon $N$ are known, the MPC problem can be solved as a single simple quadratic program (QP), for which efficient solvers are available to achieve optimal real-time solutions.
However, due to the dependence of $p$ on the
state and/or the input, unavailable before solving the optimization problem,  different techniques have been employed to handle the prediction of the scheduling parameter with successful applications, see \cite{Morato2020}.

To ensure low computational complexity, the underlying idea of LPVMPC is to solve its optimization problem as a single QP. To address the unknown values of   $p$ over $N$, the following methods have been introduced:
(1) They can be frozen at the current measured value. This leads to the so-called gain-scheduling predictive control method, e.g.,  \cite{GiAbBoGrVeMe21}.
 (2) They can be derived from the predictions of the MPC solution at the previous time step \citep{Maryam,alcala2020lpv,Karachalios2024CCTA}.
 (3) They can be approximated based on local Taylor expansions of
the scheduling function that defines $p$ as a function of the state/input  \citep{Morato23}.
 (4) They can be considered uncertain, with an uncertainty bound based on the rate of variation of $p$. Typically, robust MPC addresses such uncertainty, e.g., tube-based MPC in \cite{HanemaTubes,Abbas24}.
 The latter approach is usually the most computationally expansive due to the online construction of tubes \citep{KarachaliosICARCV2024}; however, it provides theoretical guarantees on stability and recursive feasibility.
An alternative approach was introduced in \cite{Cisneros_2016}, the qLMPC, based on solving the MPC optimization problem as a sequence of QPs based on method (2), updating the values of $p$ over $N$   iteratively until convergence.
The qLMPC performed better than the LPVMPC with a single QP at the expense of increasing the computational burden, which is still less than casting the problem as an NLP. The convergence properties have been studied in \cite{Hespe2021}. Since it solves a sequence of QPs at each MPC time interval, it resembles a sequential quadratic programming (SQP) method for solving NLPs \citep{NoceWrig06}. However, it does not directly provide insight into this connection to understand better the differences between solving NMPC as an NLP and using the LPVMPC approach, which can help evaluate and understand the suboptimality when considering the LPVMPC.

This paper considers the LPV framework for efficiently and cost-effectively solving NMPC, adopting a multi-QP approach similar to qLMPC. However, we introduce a more efficient and pragmatic alternative to qLMPC by directly leveraging LPVMPC with SQP \citep{NoceWrig06} to solve NMPC optimization problems. This can provide better insight into the differences between NLP-based NMPC and LPVMPC. Moreover, it can also incorporate other SQP tools, such as the globalization of SQP methods \citep{NoceWrig06}, leading to faster convergence with guarantees. We will consider both non-condensed and condensed forms of the optimization problem.   Moreover, we demonstrate the improvement of the solution based on the proposed approach through practical applications in simulations, showing a reduction in computational complexity, significantly faster convergence, and solutions that closely match those of the original NLP compared to state-of-the-art algorithms.

The rest of the paper is organized as follows: preliminaries and the problem setup in Sec.~\ref{sec:pre}, the proposed approach is presented in Sec.~\ref{s:propArpproach},  numerical examples in Sec.~\ref{sec:exam}, and conclusions in Sec.~\ref{sec:con}. 


\section{Preliminaries and Problem Setup}\label{sec:pre}
\subsection{LPV-based MPC and QP Formulations}
An LPV discrete-time system is represented as
\begin{subequations}\label{eq:LPV}
    \begin{align}
        x_{k+1}&=A(p_k)x_k+B(p_k)u_k,\\
        p_{k}&=\rho(x_k,u_k), 
    \end{align} 
\end{subequations}
where $x_k\in\mathbb{R}^{n}$ is the state of the system, $u_k\in\mathbb{R}^{m}$ is the input, and $p_k\in\mathcal{P}\subset\mathbb{R}^{n_\text{p}}$ is a measurable, time-varying scheduling parameter that 
belongs to a compact and convex set $\mathcal{P}$ of its admissible values for all $k\in\mathbb{Z}_+$. The system matrices 
 $A:\mathcal{P}\to\mathbb{R}^{n\times n},~B:\mathcal{P}\to\mathbb{R}^{n\times m}$   are assumed to be bounded due to the compactness of $\mathcal{P}$ with affine dpendence on $p_k$, and $\rho\in\mathcal{C}^2(\mathcal{P})$. It is assumed in this work that   the LPV representation \eqref{eq:LPV} provides  an exact embedding 
 of the dynamics of an NL system. For details about the embedding of NL dynamics in the form \eqref{eq:LPV}, see \cite{KwBoWe06} and \cite{To10}. 
 
The LPVMPC is carried out   by solving  the following optimization problem online at each  time instant $t_k\in\R_+$, where $t_k = k \cdot t_s$ and $t_s$ denotes the sample time,
%
\begin{subequations}\label{eq:opt1}
\begin{align}
    \min_{u_0, u_1,\cdots,x_0, x_1,\cdots}&\quad x_N^\top Q x_N + \sum_{k=0}^{N-1}   x_{k}^\top Qx_{k}+u_{k}^\top R u_{k} 
    \\
    \text{s.t.}&\quad x_{0}=x(t_k)\\
    &\quad x_{k+1}=A(p_{k})x_{k}+B(p_{k})u_{k}\\
    &\quad x_{k}\in\mathcal{X},~k=1,~\ldots,N\\
    &\quad u_{k}\in\mathcal{U},~k=0,~\ldots,N-1
\end{align}
\end{subequations}
where $Q\succeq 0$ and $R\succ 0$  are the quadratic cost weights, $x(t_k)$ is the current state of the system and $\mathcal{X}$, $\mathcal{U}$ 
 the state and input constraint sets are receptively represented by linear inequalities, see \eqref{eq:constraints}. Terminal conditions can be included in \eqref{eq:opt1} to ensure recursive feasibility of the optimization problem and asymptotic stability of the closed-loop system with LPVMPC, typically computed offline using LMIs, see \cite{Morato2020} for more details; however, they are not considered in this work. Since the LPV model (\ref{eq:opt1}c)  is assumed to be an exact embedding of nonlinear dynamics, by replacing $p$ with $\rho(x,u)$,  the optimization problem \eqref{eq:opt1} becomes equivalent to an NMPC. 

Two formulations are presented below to introduce the proposed approach for solving \eqref{eq:opt1}.

\subsubsection{Noncondensed Form}
By introducing  the following augmented vectors:
 \begin{equation}\label{e:augvectors}
  \bar{x}=\left[\begin{matrix}
      x_{1}\\x_{2}\\\vdots\\ x_{N}
  \end{matrix}\right],\quad \bar{u}=\left[\begin{matrix}
      u_{0}\\u_{1}\\\vdots\\ u_{N-1}
  \end{matrix}\right],\quad \bar{p}=\left[\begin{matrix}
      p_{0}\\p_{1}\\\vdots\\ p_{N-1}
  \end{matrix}\right],
 \end{equation}
we can reformulate  \eqref{eq:opt1}  in a sparse, structured optimization problem where the underlying decision variables depend on both the input and state vectors. Moreover,  define the augmented input-state decision vector as $z=[\bar{u}^\top ~~\bar{x}^\top]^\top\in\R^{N(m+n)}$. The  optimization problem  \eqref{eq:opt1} becomes 
\begin{subequations}\label{e:lpvmpcsim} 
    \begin{align}
    \min_{z}&\quad \frac{1}{2} z^\top {M} z\\
     \text{s.t.}&\quad {C}(\bar{p})z=b_{x_0},\\
     & \quad {G}^{\rm z}z\leq {h}^{\rm z},
\end{align}
\end{subequations}
where $M$ and $D$ 
are defined, respectively, as
\begin{align*}
    {M} &:=2\cdot\text{blkdiag}\left( \mathbb{I}_N\otimes R,\; \mathbb{I}_{N}\otimes Q\right),
\\
b_{x_0} &:=\begin{bmatrix}
-A^\top(p_{0})x_0  &         0&         \cdots &         0
\end{bmatrix}^\top,
\end{align*}
 and  $C(\bar{p})$, ${G}^{\rm z}$, ${h}^{\rm z}$ 
 are given in Appendix~\ref{a:NonCondMat}. The optimization problem \eqref{e:lpvmpcsim} is now a function of $\bar{p}$ and since
  $M\succeq 0$,  if $\bar{p}$ is known,  \eqref{e:lpvmpcsim} becomes a standard QP, which can be solved at each time instant.

\subsubsection{Condensed Form}
By expressing $\bar{x}$ in \eqref{e:augvectors} in terms of $x_0,\bar{u}$ the perdition model over $N$ is given as
\begin{equation}\label{e:AugLPVSys}
  \bar{x}=\bar{A}(\bar{p})x_0+\bar{B}(\bar{p})\bar{u},
\end{equation}
 where  $\bar{A},~\bar{B}$ are  augmented matrices as 
 given in Appendix~\ref{a:NonCondMat}.
Consequently, the cost function (\ref{eq:opt1}a) can be written solely in $\bar{u}$. By collecting all the input and state constraints in one inequality constraint of $\bar{u}$, the constrained optimization problem \eqref{eq:opt1} can be reformulated as 
\begin{subequations}\label{e:LPVMPCQP}
    \begin{align}
        \min_{\bar{u}} &\quad \frac{1}{2} \bar{u}^\top H(\bar{p})\bar{u}+ x_0^\top g(\bar{p})\bar{u} \\
        \text{s.t. }&\quad    {G}^{\rm u}(\bar{p})\bar{u}\leq {h}^{\rm u}(\bar{p}),
    \end{align}
\end{subequations}
where
\[
H(\bar{p})  \!=\!\begin{bmatrix}
\mathbb{I}_{Nm}\\ \bar{B}(\bar{p})
\end{bmatrix}^\top \!\! M\!\begin{bmatrix}
\mathbb{I}_{Nm}\\ \bar{B}(\bar{p})
\end{bmatrix},\;\;
g(\bar{p})\!=\!\begin{bmatrix}
0_{Nm}\\ \bar{A}(\bar{p})
\end{bmatrix}^\top \!\!M\!\begin{bmatrix}
0_{Nm}\\ \bar{B}(\bar{p})
\end{bmatrix}
\]
and ${G}^{\rm u}(\bar{p})$, ${h}^{\rm u}(\bar{p})$ are given in Appendix~\ref{a:NonCondMat}. With $H(\bar{p})\succeq 0$, $\forall \bar{p}$,    if $\bar{p}$ is known,  \eqref{e:LPVMPCQP} becomes  a standard QP.

Generally, the noncondensed form is numerically superior if the underlying solver can exploit sparsity; otherwise, the condensed form is preferable. The formulations \eqref{e:lpvmpcsim}  and \eqref{e:LPVMPCQP} are used in the next section for introducing the proposed approach.

%


 
\subsection{Basic SQP for solving  NLP} 
 
Next, the basic SQP method  \citep{NoceWrig06}  is reviewed. It is widely used for solving NLP problems in optimal control and NMPC.  Under certain regularity conditions, SQP guarantees local convergence to a solution of the NLP problem. Next, we illustrate the  idea of the SQP approach using    the noncondensed form \eqref{e:lpvmpcsim}\footnote{Similar analysis can be performed using the condensed form \eqref{e:LPVMPCQP}.}, which, since $p=\rho(x,u)$, it can be reformulated in terms of  $z$ as follows
\begin{subequations}    \label{e:NLP}
    \begin{align}
    \min_{z}&\quad \frac{1}{2} z^\top {M} z\\
     \text{s.t.}&\quad {C}(\bar{p}(z))z=b_{x_0},\\
     & \quad {G}^{\rm z}z\leq {h}^{\rm z}.
     \end{align}
\end{subequations}
Therefore, \eqref{e:NLP} is an NLP and can be solved using SQP.
The principle idea is to successively linearize the   Karush-Kuhn-Tucker (KKT) optimality conditions of \eqref{e:NLP} to construct a sequence of QP subproblems. In the following, the dependence of  $\bar{p}(z)$ is omitted for simplicity of notation. The Lagrangian of the   NLP \eqref{e:NLP} is given as
\begin{equation}\label{e:Lag-NLP}
    \La(z,\lambda)=\frac{1}{2} z^\top {M} z+\lambda^\top({C}(z)z-{b}_{x_0})
          +\mu^\top({G}^{\rm z}z- {h}^{\rm z}),
\end{equation}
where $\lambda$ and $\mu$ are the Lagrange multipliers for the equality and inequality constraints, respectively. The SQP method then iteratively solves \eqref{e:NLP},  yielding the updates 
\begin{equation}\label{e:itupdate-NLP}
    z^{l+1}=z^{l}+d^l,\quad\lambda^{l+1}=\delta^l,\quad \mu^{l+1}=\eta^l,
    \end{equation}
    where the current iterate 
    $(d^l,\delta^l,\eta^l)$ is the solution of the following QP subproblem 
    \begin{subequations}    \label{e:NLP2QP}
    \begin{align}
    \min_{d^l}&\quad\frac{1}{2}{d^l}^\top \nabla^2_{z}\La(z^l,\lambda^l,\mu^l) d^l+ \left(MZ^l\right)^\top d^l
     \\
     \text{s.t.}&\quad \frac{\partial [{C}(z^l)z^l]}{\partial z}d^l=b_{x_0}-{C}(z^l)z^l,\\
     &\quad {G}^{\rm z}d^l\leq {h}^{\rm z}-{G}^{\rm z}z^l,
\end{align}
\end{subequations}
where $ \nabla^2_{z}\La(z^l,\lambda^l,\mu^l)$ is the Hessian
of the Lagrangian \eqref{e:Lag-NLP} and  $\frac{\partial [{C}(z^l)z^l]}{\partial z}$ represents  the  Jacobian of the equality constraint  
(\ref{e:NLP}b), evaluated at the iterate ($z^l,\lambda^l,mu^l$). Under certain regularity assumptions, including the continuity and differentiability of the optimization problem and the satisfaction of constraint qualifications, SQP is guaranteed to converge to the optimal local solution $(z^\ast,\lambda^\ast,\mu^\ast)$, provided that the initialization is sufficiently close to the solution, see \cite{NoceWrig06} for more details.  However, the computational complexity of computing the Hessian and Jacobian at every iteration is one of the main challenges of SQP.  To mitigate this, inexact formulations of these terms are commonly considered to reduce this complexity as shown in \cite{Bock2007a}, see also \cite{Bock2007a,RaMaDi17chapter8}. However, this often leads to suboptimal solutions. 

In this work, based on the LPVMPC formulation, we present a systematic approach that can be seen as solving SQP using inexact Hessians and Jacobians, which can be easily updated using the scheduling parameter, enabling inexpensive construction of the QP subproblems and leading to suboptimal solutions closer to the optimal ones.

\section{Proposed Approach}\label{s:propArpproach}
In this section, we propose an approach based on the LPV  formulation to solve the NLP associated with \eqref{eq:opt1} using approximate SQP methods,  where inexact Hessians and Jacobians depend on the scheduling parameter according to the LPV embedding used to represent the nonlinear model. We apply this approach to both the noncondensed and the condensed forms.

Consider  the noncondensed form \eqref{e:NLP} and  introduce the  following inexact Hessian and Jacobian for \eqref{e:NLP2QP} 
\begin{equation}\label{e:approxHessGrad-noncond}
 \nabla^2_{z}\La(z^l,\lambda^l,\mu^l)\approx {M},\quad \frac{\partial [{C}(\bar{p}(z^l))z^l]}{\partial z}\approx {C}(\bar{p}^l),
\end{equation}
where $\bar{p}^l=\bar{p}(z^l)$.
This  leads to the following approximate QP subproblem
\begin{subequations}    \label{e:QP-variant-2}
    \begin{align}
    \min_{d^l}&\quad   \frac{1}{2}   {d^l}^\top  {M} d^l+({M}z^l)^\top d^l
     \\
     \text{s.t.}&\quad  {C}(\bar{p}^l)d^l=b_{x_0}-{C}(\bar{p}^l)z^l,\\
    &\quad  {G}^{\rm z}d^l\leq {h}^{\rm z}-{G}^{\rm z}z^l,
\end{align}
\end{subequations}
which should be solved iteratively given $z^l$ and $\bar{p}^l$,  updating $z^{l+1}=d^l + z^l$ and $\bar{p}^{l+1}$ accordingly  until $d^l$ converges to zero. It is worth mentioning that if such an iterative QP converges, it will not exactly reach the solution of \eqref{e:NLP} but rather a suboptimal one. Conditions similar to those used for analyzing the convergence properties of NLP algorithms with inexact Hessians and Jacobians, as discussed in \cite{Bock2007a}, can be applied to show convergence proprieties of \eqref{e:QP-variant-2}, see also \cite{Hespe2021}. However, this paper focuses on introducing the approach, while the convergence properties will be addressed in a future study.

Note that adding the constant term ${z^l}^\top {M}z^l$ to the cost function (\ref{e:QP-variant-2}a) 
does not affect the solution; however,  it allows the expression to be written as ${z^{l+1}}^\top {M}z^{l+1}$. Then, by substituting $d^l+z^l$ with $z^{l+1}$ in (\ref{e:QP-variant-2}b,c),  we recover the   LPVMPC nonconduced form in \eqref{e:lpvmpcsim}, where the QP  being directly scheduled with $\bar{p}$.  This provides insight into the proposed inexact Hessian and Jacobian in \eqref{e:approxHessGrad-noncond}, enabling a systematic and cost-effective construction and update of the QP subproblem without explicitly forming these inexact components.
It is not just a simple approximation; it also offers meaningful insight into the control problem, as it allows the inexact Hessian and Jacobian to adapt to the behavior of the controlled system through the choice of scheduling parameter.

Next, we consider the condensed form in \eqref{e:LPVMPCQP} in its NLP formulation by directly rewriting   it  in terms of $\bar{u}$, since $\bar{p}=\bar{p}(\bar{u})$,     as follows: 
\begin{subequations}\label{e:NLPseq}
    \begin{align}
        \min_{\bar{u}}&\quad \frac{1}{2} \bar{u}^\top H(\bar{u})\bar{u}+ x_0^\top g(\bar{u})\bar{u}, \\
        \text{s.t. }&\quad    {G}^{\rm u}(\bar{u})\bar{u}\leq {h}^{\rm u}(\bar{u}).
    \end{align}
\end{subequations}
Let $a(\bar{u})= \frac{1}{2} \bar{u}^\top H(\bar{u})\bar{u}+ x_0^\top g(\bar{u})\bar{u}$, then the Lagrangian of \eqref{e:NLPseq} is given by
\begin{equation}\label{e:Lagrangian_NLPseq}
\La(\bar{u},\mu)= a(\bar{u})+\mu^\top \left( {G}^{\rm u}(\bar{u})\bar{u} - {h}^{\rm u}(\bar{u})\right).
\end{equation}
The  SQP method for solving the NLP \eqref{e:NLPseq} results in an iteration update by solving the following QP subproblem
 \begin{subequations}\label{e:NLPseqtoQP}
    \begin{align}
        \min_{d^l} &\quad\frac{1}{2}{d^l}^\top\nabla_{\bar{u}}^2
        \La(\bar{u}^l,\mu^l) d^l +\nabla_{\bar{u}} a(\bar{u}^l)^\top d^l
          \\
        \text{s.t.}&\quad \frac{\partial[{G}^{\rm u}(\bar{u}^l)\bar{u}^l-{h}^{\rm u}(\bar{u}^l)]}{\partial \bar{u}}  d^l \leq {h}^{\rm u}(\bar{u}^l)- {G}^{\rm u}(\bar{u}^l)\bar{u}^l.
    \end{align}
\end{subequations}
Under a regularity assumption for the solution of \eqref{e:NLPseq}  and with exact Hessian, gradient, and Jacobian, it is guaranteed that solving \eqref{e:NLPseqtoQP} iteratively will converge to the optimal local solution of \eqref{e:NLPseq}, provided a reasonable initialization close to the solution. However, the exact Hessian of the Lagrangian, the gradient of the cost function, and the Jacobian of the inequality constraints are quite complex,  highlighting the difficulty of solving the NLP \eqref{e:NLPseq} based on the QP subproblem in  \ \ref {e:NLPseqtoQP}.  To avoid such complexity, we propose the following approximations. 
\begin{subequations}\label{e:approxCondens}
\begin{align}
\nabla_{\bar{u}}^2\La(\bar{u}^l,\mu^l)&\approx  H(\bar{p}^l),\\
 \nabla_{\bar{u}}f(\bar{u}^l)&\approx H(\bar{p}^l)\bar{u}^l+g^\top(\bar{p}^l)x_0,\\
 \frac{\partial[{G}^{\rm u}(\bar{u}^l)\bar{u}^l-{h}^{\rm u}(\bar{u}^l)]}{\partial \bar{u}}&\approx {G}^{\rm u}(\bar{p}^l),
\end{align}
\end{subequations}
where the dependence in the right-hand side is now written in terms of $\bar{p}$, since $\bar{p}=\bar{p}(\bar{u})$, and the explicit dependence on $\bar{u}$  is dropped.
These approximations yield the following condensed QP subproblem
 \begin{subequations}\label{e:seq-var1}
    \begin{align}
        \min_{d^l}&\quad \frac{1}{2} {d^l}^\top H(\bar{p}^l) d^l+ \left(H(\bar{u}^l)\bar{p}^l+g^\top (\bar{p}^l)x_0\right)^\top d^l
          \\
        \text{s.t.}&\quad
{G}^{\rm u}(\bar{p}^l)d^l  \leq {h}^{\rm u}(\bar{p}^l)- {G}^{\rm u}(\bar{p}^l) \bar{u}^l. 
    \end{align}
\end{subequations}
As demonstrated in the noncondensed form, it is straightforward to show how \eqref{e:LPVMPCQP} can be recovered from \eqref{e:seq-var1}, leading to the intuition behind the approximations considered in \eqref{e:approxCondens}.
Considering \eqref{e:seq-var1} iteratively for solving \eqref{e:NLPseq}, convergence to a suboptimal solution holds once $d^l\approx 0$. The convergence properties can also be analyzed using the approach of \cite{Bock2007a}.

It is important to note that iteratively solving the QP subproblem in   \eqref{e:QP-variant-2}/\eqref{e:seq-var1}, which, if it converges,  yields a solution that differs from that of \eqref{e:NLP}/\eqref{e:NLPseq} due to the approximations employed. However, it turns out that, through several numerical applications, the suboptimal solutions based on \eqref{e:QP-variant-2}/\eqref{e:seq-var1} usually lead to solutions very close to the optimal ones, see Section~ \ref{sec:exam}. 


\begin{remark}
It is worth mentioning that the qLMPC approach of \cite{Cisneros_2016} is based on iteratively solving QPs, with updates of $z^{l+1} / \bar{u}^{l+1}$ according to the condensed form. In each iteration, $\bar{p}^{l}$ is used to compute the current QP subproblem. However, the convergence of the qLMPC is achieved once the difference between $\bar{p}^{l+1}$ and $\bar{p}^{l}$ is smaller than a given stopping criterion. The convergence properties of this scheme have been studied in \cite{Hespe2021}. The key difference between the proposed approach here and the qLMPC is that solving the QP subproblem in our approach is based on the increment $d^l$, as proposed in \eqref{e:QP-variant-2} and \eqref{e:seq-var1} leading to a numerically robust scheme and can yield better suboptimal solutions, as demonstrated in the numerical simulations presented in the next section. Furthermore, it provides better insight into the problem, its connection to SQP, and a meaningful understanding of the difference between solving the NLP and considering the LPVMPC approach with the approximate SQP.
\end{remark}



\section{NUMERICAL EXAMPLES}\label{sec:exam}

This section tests the proposed approach, referred to as  LPVMPC-SQP, against established benchmarks. 
The simulations are performed on a Dell Latitude 5590 laptop with an Intel(R) Core(TM) i7-8650U CPU and 16 GB of RAM. They are implemented in MATLAB. The LPVMPC-SQP approach is implemented using the {\tt quadprog} command based on the {\tt active-set} algorithm to solve the QP subproblem; therefore, execution follows the condensed form \eqref{e:seq-var1}.
\subsection{Example 1 - The forced Van der Pol oscillator}

We start with the forced Van der Pol oscillator, a well-studied and challenging benchmark \citep{VaHo14}.
 The continuous-time nonlinear ordinary differential equation of the system  is given by
\begin{equation}\label{e:vdp}
    \ddot{y}(t)=(1-y^2(t))\dot{y}(t)-y(t)+u(t),
\end{equation}
where $y$ is the output of the system, subject to the constraints $\dot{y}\geq -0.25$, $u$ is the input constrained by $|u|\leq 1$. Let $x=[y\;\;\dot{y}]^\top$ be the state vector. To implement MPC, the model is discretized using the classical fourth-order Rung-Kutta method with a sampling time of $t_\text{s}=0.5$s. For the proposed schemes, the LPV embedding is derived as in \eqref{eq:LPV}, where both system matrices $A$ and $B$ are parameter-dependent with an affine dependence on 11 scheduling parameters. The components of the scheduling vector depend on both $x$ and $u$. The explicit expressions for them are omitted here due to space limitations.
For comparison, we consider NMPC solved using the NLP solver {\tt IPOPT} \citep{WaBi06}, which is implemented in {\tt CasADi} \citep{andersson2018casadi} and referred to here as NMPC-IPOPT.

\begin{figure}
    \centering
    \includegraphics[scale=0.5]{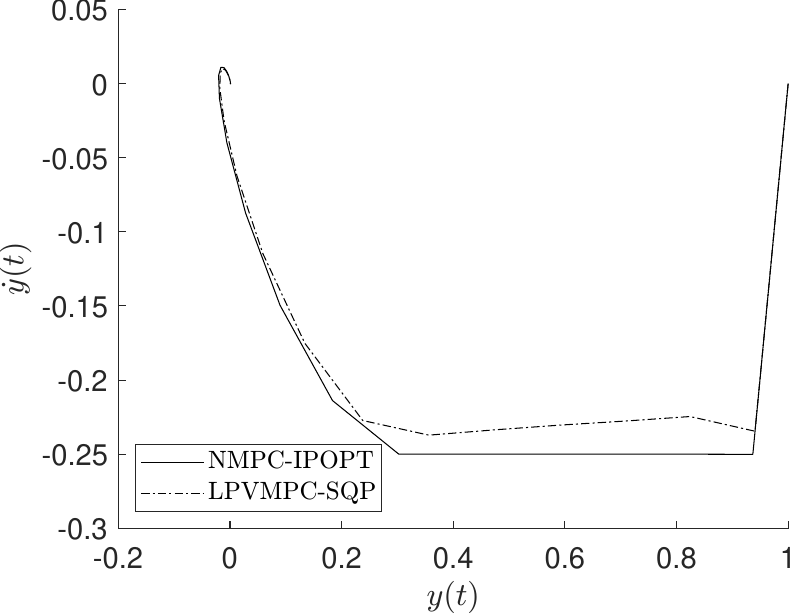} 
\caption{Phase plot  of the closed-loop state trajectories of the Van der Pol oscillator.}
    \label{fig:phase}
\end{figure}
\begin{figure}
    \centering
\includegraphics[scale=0.5]{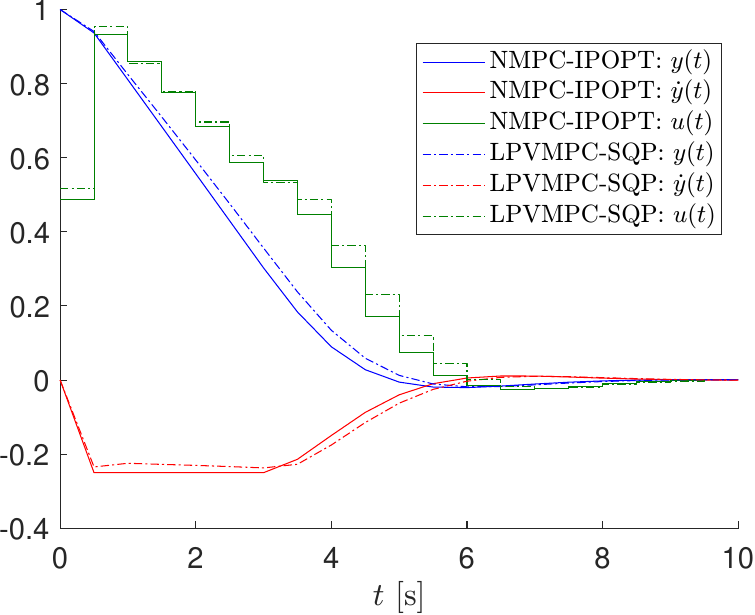}
\caption{The closed-loop state and input trajectories  of  the Van der Pol oscillator.}
    \label{fig:traj}
\end{figure}
\begin{figure}
    \centering
    \includegraphics[scale=0.5]{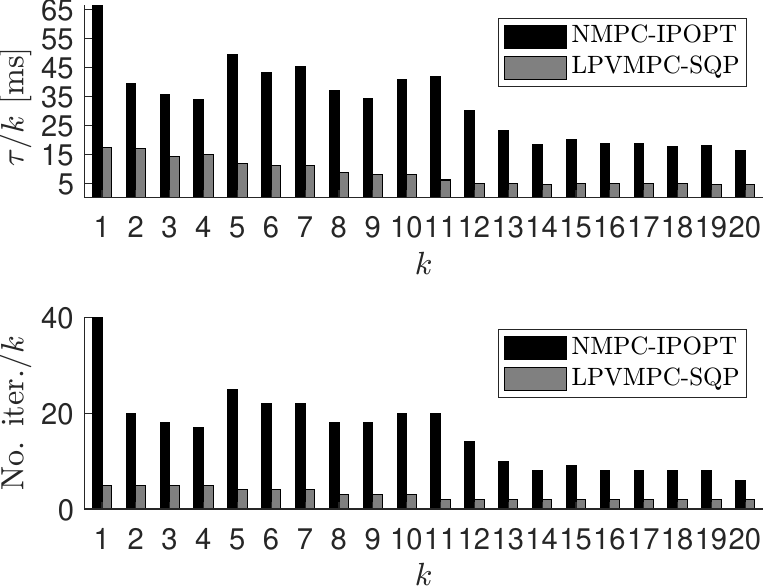} 
\caption{Solver execution time and iterations per MPC  step for the considered approaches applied to the Van der Pol oscillator averaged over 100 simulation runs.}
    \label{fig:solperform}
\end{figure}

The control objective is to steer the
initial conditions $x_0=[1~~0]^\top$   to the origin quickly while respecting the given input/state constraints. 
We consider the implementation of both MPC approaches on the continuous-time system \eqref{e:vdp}, where the decision variables in both schemes  are initialized to zero at each MPC execution step.
The resulting closed-loop state trajectories are shown in the phase plot in Fig.~\ref{fig:phase}, and the time trajectories of the states and input are given in Fig.~\ref{fig:traj}. To evaluate optimality, the total cost of the optimal control problem, based on the closed-loop trajectories, is 6.9612 for the NMPC-IPOPT approach and 7.3684 for the proposed LPVMPC-SQP approach. This indicates that  LPVMPC-SQP  results in an optimality loss of approximately 5.85\% compared to the optimal cost achieved by  NMPC-IPOPT,  which is reasonably satisfactory given the approximations of the Hessian and Jacobian used in the formulation of LPVMPC-SQP. Given this slight difference, the system's performance using LPVMPC-SQP closely matches the optimal performance achieved by the NMPC-IPOPT approach, as illustrated in Figs.~\ref{fig:phase} and \ref{fig:traj}. Interestingly, the computational burden of the LPVMPC-SQP, in terms of execution time $\tau$ per MPC execution step $k$ and the number of iterations per MPC step, is on average (over 100 simulation runs) reduced by approximately 73.53\% in execution time and 80.88\% in the number of iterations compared to the NMPC-IPOPT, as shown in Fig.~\ref{fig:solperform}. This significant reduction in computational cost is achieved at the expense of only an optimality loss of  5.85\%, demonstrating the efficiency of the proposed method.

\subsection{Example 2 - The dynamic unicycle}
The second example concerns a dynamic unicycle system borrowed from \cite{Hespe2021}, where a similar case study was conducted. We provide further comparisons with the proposed methods. The unicycle model represents a wheeled robot with simplified dynamics,  subject to non-holonomic constraints. The nonlinear dynamics, based on its kinematics,  can be described as:
\begin{equation}\label{eq:unicycle}
    \begin{aligned}
        \dot{s}(t)&=\upsilon(t)\cos(\phi(t)),\\
        \dot{q}(t)&=\upsilon(t)\sin(\phi(t)),\\
        \dot{\upsilon}(t)&=F(t),\\
        \dot{\phi}(t)&=\omega(t),\\
        \dot{\omega}(t)&=r(t).
    \end{aligned}, 
\end{equation}
where $s,q$ denotes the robot's positions in the global coordinate system, $\phi$ is the orientation, $\upsilon$ is the forward velocity, and $\omega$ is the rotational velocity.
The model \eqref{eq:unicycle} is discretized using the  forward Euler method with sampling time $t_s=0.1$s, and the discrete-time 
state $x_k:=\left[s_k,q_k,\upsilon_k,\phi_k,\omega_k\right]^\top$, input $u_k=\left[F_k,r_k\right]^\top$, and  scheduling parameter vector $p_k=\left[\cos(\phi_k),\sin(\phi_k)\right]^\top$. The 
LPV embedding of the nonlinear dynamics \eqref{eq:unicycle}, expressed in the  
form of \eqref{eq:LPV}, is given as: 
\begin{equation*}
    A(p_k)=\left[\begin{matrix}
        1 & 0 & t_\text{s}\cos(\phi_k) & 0 & 0\\
        0 & 1 & t_\text{s}\sin(\phi_k) & 0 & 0\\
        0 & 0 & 1 & 0 & 0\\
        0 & 0 & 0 & 1 & t_\text{s}\\
        0 & 0 & 0 & 0 & 1
    \end{matrix}\right],~B=\left[\begin{matrix}
        0 & 0\\
        0 & 0\\
        t_\text{s} & 0\\
        0 & 0\\
        0 & t_\text{s}
    \end{matrix}\right].
\end{equation*}

\begin{figure}[!h]
    \centering
  \hspace{0mm}  \includegraphics[width=1\linewidth]{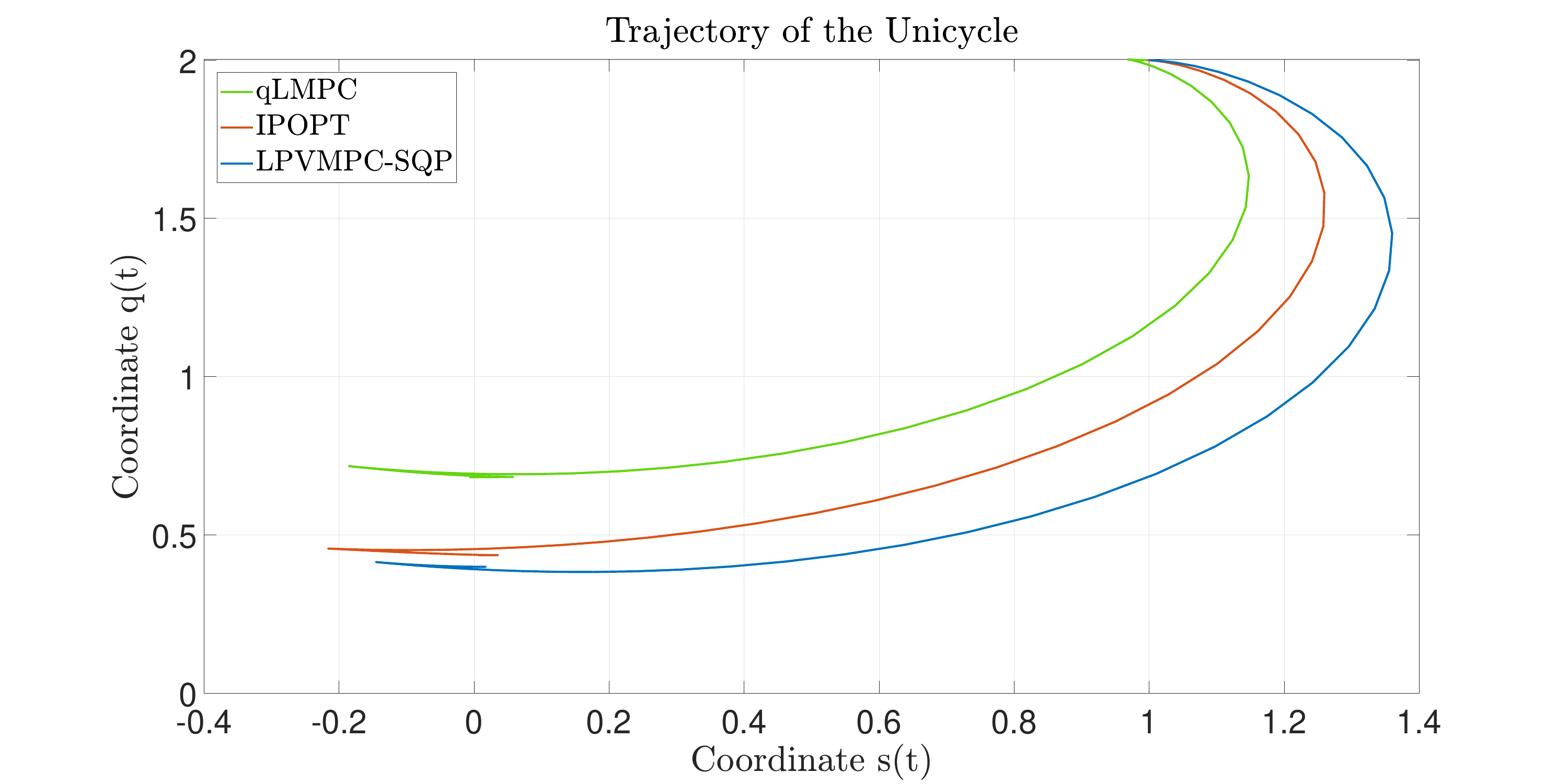}
    \caption{Closed-loop trajectories of the regulation problem for the unicycle.}
    \label{fig:fig1_Unicycle}
\end{figure}

The unicycle's initial condition is $x_0=\left[1,2,0,\pi,0\right]^\top$, and the task is to regulate the state within 10 s 
to the origin as in Fig.~\ref{fig:fig1_Unicycle}. 
For comparison, the qLMPC approach of \cite{Cisneros_2016} 
and the NMPC, with its NLP, is solved 
using {\tt IPOPT} implemented in {\tt CasADi}.
For all approaches, the  MPC setting is as follows: prediction horizon  $N=20$, and the quadratic weighted matrices   $Q=\texttt{diag}(1,1,0.1,1,0.1)$ and $R=\texttt{diag}(1,1)$. The condensed form is used for all MPC schemes.

The obtained closed-loop trajectories for the different methods are depicted on the phase plane with global coordinates in  Fig.~\ref{fig:fig1_Unicycle}. None of the tested methods succeeded in steering the unicycle at the origin, 
due to the non-holonomic constraints. However, 
the proposed method LPVMPC-SQP  converged to almost the optimal point, comparable to  NMPC-IPOPT 
as depicted in Fig.~\ref{fig:fig1_Unicycle}, outperforming qLMPC.  
The solver timing in Table~\ref{tab:Compuation_time2} indicates that the proposed LPVMPC-SQP method and qLMPC provide similar computation performance, in contrast with the $\sim 10$ times slower NMPC-IPOPT. 


       \begin{table}[!h]
    \centering
    \caption{Comparison of computation times for the unicycle example. 
    The average execution time is denoted $\bar\tau$. 
    }
    \label{tab:Compuation_time2}
    \setlength{\tabcolsep}{2pt}
    \begin{tabular}{ l |  c | c | c } 
        \toprule
        \rowcolor{mydarkgray} \textbf{Controller}  & \textbf{$\bar{\tau}$}  & \textbf{$\sigma_{\tau}$}  & \textbf{$\tau_{max}$} \\ \hline
                                    IPOPT &  0.1008   &  0.1577  &  1.6384\\
        \rowcolor{mygray}   LPVMPC-SQP &    0.0037  &  \textbf{0.0013} &  \textbf{0.0080}\\
         qLMPC &     \textbf{0.0035}  &  0.0016  &  0.0082\\ 
        \end{tabular}
        \end{table}
\subsection{Example 3 - Reference tracking in Autonomous Driving}
The third example concerns modeling the car dynamics with the dynamic bicycle model described in \cite{Maryam} for solving a reference tracking problem in the context of autonomous driving. In this example, we only consider the proposed approach. 
The continuous time differential equations describing the motion of a vehicle are presented as
\begin{equation}\label{eq:dynamic_model}
\begin{aligned}
    \dot{X}(t) &= \upsilon(t)\cos{\psi(t)}-\nu(t)\sin{\psi(t)}, \\
    \dot{Y}(t) &= \upsilon(t)\sin{\psi(t)}+\nu(t)\cos{\psi(t)},  \\
    \dot{\upsilon}(t) &= \omega(t) \nu(t) + a(t), \\
    \dot{\nu}(t) &= -\omega(t)\upsilon(t) + \frac{2}{m}( F_{\rm yf}(t) \cos{\delta(t)} + F_{\rm yr}(t)), \\
    \dot{\psi}(t) &= \omega(t), \\
    \dot{\omega}(t) &= \frac{2}{I_{\rm z}} ( l_{\rm f} F_{\rm yf}(t) - l_{\rm r}  F_{\rm yr}(t)),
\end{aligned}
\end{equation}
where $X$, $Y$, $\upsilon$, $\nu$, $\psi$ and $\omega$ denote the global $X$ axis coordinate of the center of gravity (GoG), the global $Y$ axis coordinate of the CoG, the longitudinal speed in body frame, the lateral speed in body frame, the vehicle yaw angle and the yaw angle rate, respectively. The control inputs of the system are the longitudinal acceleration $a$ and the steering angle $\delta$. 
The vehicle moment of inertia and mass are denoted by $I_{\rm z}$ and $m$, respectively. The lateral forces acting on the front and rear tires are denoted as $F_{\rm yf}$ and $F_{\rm yr}$, respectively, and calculated linearly as $F_{\rm yf} = C_{\alpha \rm f}  \alpha_{\rm f}$, $F_{\rm yr} = C_{\alpha \rm r}  \alpha_{\rm r}.$
The parameters $C_{\alpha \rm f}$ and $C_{\alpha \rm r}$ represent the cornering stiffness of the front and rear tires, respectively. The slip angle of the front tire is $\alpha_{\rm f}$ and is calculated as $\alpha_{\rm f} = \delta - (\nu + l_{\rm f} \omega)/\upsilon$. The rear tire slip angle is $\alpha_{\rm r}$ and is calculated as $\alpha_{\rm r} = (l_{\rm r} \omega -\nu )/\upsilon$. The above nonlinear system can be embedded in the LPV as
\begin{equation}\label{eq:dis_model}
\left\{\begin{aligned}
    \dot{z}(t)&=A_c(p(t))z(t)+B_c(p(t))u(t),\\
    p(t)&=(\upsilon(t),\nu(t),\delta(t),\psi(t)),~t\geq 0.
    \end{aligned}\right.
\end{equation}
The car parameters and the LPV embedding is the same as in \cite{Maryam} where our aim in this study is to solve the optimal control problem with the newly proposed LPVMPC-SQP method.  the  MPC setting is as follows: prediction horizon  $N=15$, sampling time $t_s=0.05$ (s) and the quadratic weighted matrices   $Q=\texttt{diag}(1,1,1,1,1,1)$ and $R=\texttt{diag}(0.1,0.1)$, and the condensed form is used for the MPC.

\begin{figure}
    \centering
    \includegraphics[width=0.9\linewidth]{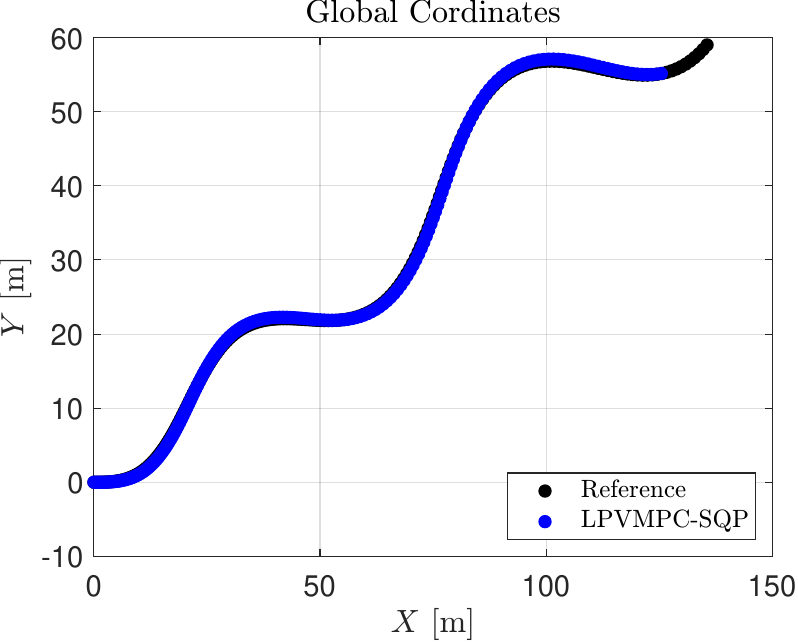}
    \caption{Phase plot of the reference tracking control problem. 
    }
    \label{fig:phasetracking}
\end{figure}

\begin{figure}
    \centering
    \includegraphics[width=0.9\linewidth]{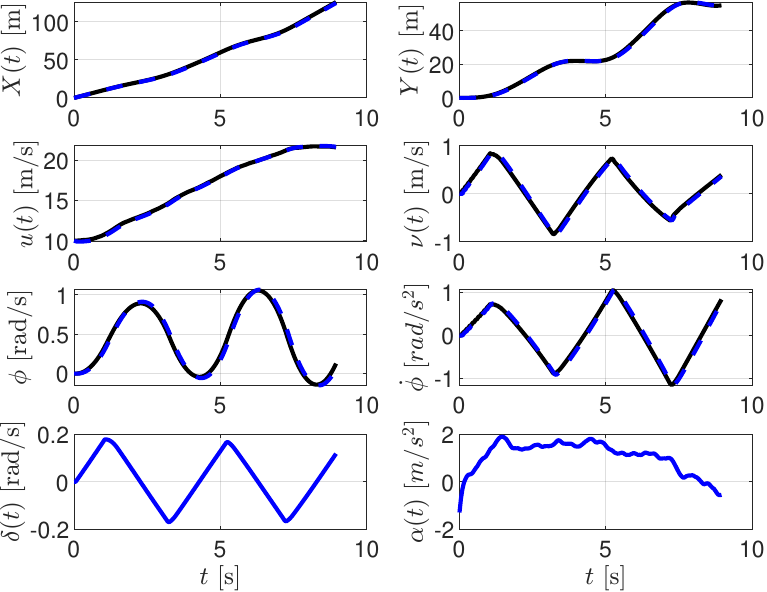}
    \caption{The state and  control trajectories of the reference tracking problem. 
    }
    \label{fig:trajtracking}
\end{figure}

The reference tracking result in the phase plot is shown in Fig.\ref{fig:phasetracking}, demonstrating high-performance tracking of the vehicle to the reference trajectory with almost negligible tracking error using the proposed LPVMPC-SQP method. The corresponding state and input trajectories are shown in Fig.\ref{fig:trajtracking}. The results obtained with LPVMPC-SQP are achieved with efficient computational time, as depicted in Table~\ref{tab:Compuation_time2}.

\begin{table}[!h]
    \centering
    \caption{Performance and computation time for the car example. 
    }
    \label{tab:Compuation_time2}
    \setlength{\tabcolsep}{2pt}
    \begin{tabular}{ l | c | c | c | c } 
        \toprule
        \rowcolor{mydarkgray} \textbf{Controller} & \textbf{Cost}-$\sum J_k$  & \textbf{$\bar{\tau}$}  & \textbf{$\sigma_{\tau}$}  & \textbf{$\tau_{max}$} \\ \hline
                                      
           LPVMPC-SQP &  \textbf{121.76}  &  \textbf{0.0107}  &  \textbf{0.0098} &  \textbf{0.1209} 
        \end{tabular}
        \end{table}
        
\section{Conclusion}\label{sec:con}

In this work, we proposed an approach for obtaining an efficient approximate solution to NMPC by combining LPVMPC with SQP. The approach leads to an associated inexact Hessian and Jacobian for the optimization problem, which depends on the scheduling parameter in the LPV embedding of the nonlinear MPC model. This enables a cost-effective update of these components throughout the SQP iterations. Based on several simulations with different benchmarks and comparisons with state-of-the-art approaches, LPVMPC-SQP demonstrates numerical advantages over NLP methods, which require costly updates of the exact Jacobian and Hessian leading to suboptimal solutions close to the optimal ones.



The approach can be extended to different SQP tools, such as globalization techniques for SQP methods, leading to even faster convergence with guarantees. Finally, further analysis is required to establish sufficient conditions for proving convergence, providing a theoretical justification for the proposed method.

\bibliography{ifacconf}             

\begin{thebibliography}{21}
\providecommand{\natexlab}[1]{#1}
\providecommand{\url}[1]{\texttt{#1}}
\providecommand{\urlprefix}{URL }
\expandafter\ifx\csname urlstyle\endcsname\relax
  \providecommand{\doi}[1]{doi:\discretionary{}{}{}#1}\else
  \providecommand{\doi}{doi:\discretionary{}{}{}\begingroup \urlstyle{rm}\Url}\fi

\bibitem[{Abbas(2024)}]{Abbas24}
Abbas, H.S. (2024).
\newblock Linear parameter-varying model predictive control for nonlinear systems using general polytopic tubes.
\newblock \emph{Automatica}, 160, 111432.

\bibitem[{Alcal{\'a} et~al.(2020{\natexlab{a}})Alcal{\'a}, Puig, and Quevedo}]{alcala2020lpv}
Alcal{\'a}, E., Puig, V., and Quevedo, J. (2020{\natexlab{a}}).
\newblock {L}{P}{V}-{M}{P} planning for autonomous racing vehicles considering obstacles.
\newblock \emph{Robotics and Autonomous Systems}, 124, 103392.

\bibitem[{Alcal{\'a} et~al.(2020{\natexlab{b}})Alcal{\'a}, Puig, Quevedo, and Rosolia}]{alcala2020autonomous}
Alcal{\'a}, E., Puig, V., Quevedo, J., and Rosolia, U. (2020{\natexlab{b}}).
\newblock {Autonomous racing using linear parameter varying-model predictive control (LPV-MPC)}.
\newblock \emph{Control Engineering Practice}, 95, 104270.

\bibitem[{Andersson et~al.(2018)Andersson, Gillis, Horn, Rawlings, and Diehl}]{andersson2018casadi}
Andersson, J.A.E., Gillis, J., Horn, G., Rawlings, J.B., and Diehl, M. (2018).
\newblock {CasADi}: A software framework for nonlinear optimization and optimal control.
\newblock \emph{Mathematical Programming Computation}, 11(1), 1--36.

\bibitem[{Bock et~al.(2007)Bock, Diehl, Kostina, and Schl\"oder}]{Bock2007a}
Bock, H.G., Diehl, M., Kostina, E.A., and Schl\"oder, J.P. (2007).
\newblock Constrained optimal feedback control of systems governed by large differential algebraic equations.
\newblock In L.T. Biegler, O.~Ghattas, M.~Heinkenschloss, D.E. Keyes, and B.~van Bloemen~Waanders (eds.), \emph{Real-Time PDE-Constrained Optimization}, 3--24. SIAM.

\bibitem[{Cisneros et~al.(2016)Cisneros, Voss, and Werner}]{Cisneros_2016}
Cisneros, P.S.G., Voss, S., and Werner, H. (2016).
\newblock Efficient nonlinear model predictive control via quasi-{L}{P}{V} representation.
\newblock In \emph{2016 IEEE 55th Conference on Decision and Control (CDC)}, 3216--3221.

\bibitem[{Gidon et~al.(2021)Gidon, Abbas, Bonzanini, Graves, Velni, and Mesbah}]{GiAbBoGrVeMe21}
Gidon, D., Abbas, H.S., Bonzanini, A.D., Graves, D.B., Velni, J.M., and Mesbah, A. (2021).
\newblock Data-driven lpv model predictive control of a cold atmospheric plasma jet for biomaterials processing.
\newblock \emph{Control Engineering Practice}, 109, 104725.

\bibitem[{Hanema et~al.(2016)Hanema, Tóth, and Lazar}]{HanemaTubes}
Hanema, J., Tóth, R., and Lazar, M. (2016).
\newblock Tube-based anticipative model predictive control for linear parameter-varying systems.
\newblock In \emph{2016 IEEE 55th Conference on Decision and Control (CDC)}, 1458--1463.

\bibitem[{Hespe and Werner(2021)}]{Hespe2021}
Hespe, C. and Werner, H. (2021).
\newblock Convergence properties of fast quasi-{LPV} model predictive control.
\newblock In \emph{2021 60th IEEE Conference on Decision and Control (CDC)}, 3869--3874.

\bibitem[{Karachalios and Abbas(2024)}]{Karachalios2024CCTA}
Karachalios, D.S. and Abbas, H.S. (2024).
\newblock Parameter refinement of a ballbot and predictive control for reference tracking with linear parameter-varying embedding.
\newblock In \emph{Proc. IEEE Conf. Control and Communication Technologies (CCTA)}, 688--693. Newcastle, UK.

\bibitem[{Karachalios et~al.(2024)Karachalios, Nezami, Schildbach, and Abbas}]{KarachaliosICARCV2024}
Karachalios, D.S., Nezami, M., Schildbach, G., and Abbas, H.S. (2024).
\newblock Error bounds in nonlinear model predictive control with linear differential inclusions of parametric-varying embeddings.
\newblock In \emph{Proc. the 18th International Conference on Control, Automation, Robotics and Vision (ICARCV 2024)}, 866--871.

\bibitem[{Kwiatkowski et~al.(2006)Kwiatkowski, Boll, and Werner}]{KwBoWe06}
Kwiatkowski, A., Boll, M., and Werner, H. (2006).
\newblock Automated generation and assessment of affine {LPV} models.
\newblock In \emph{Proc. 45th IEEE Conference on Decision and Control}, 6690--6695. San Diego, CA, USA.

\bibitem[{Morato et~al.(2020)Morato, Normey-Rico, and Sename}]{Morato2020}
Morato, M.M., Normey-Rico, J.E., and Sename, O. (2020).
\newblock Model predictive control design for linear parameter varying systems: A survey.
\newblock \emph{Annual Reviews in Control}, 49, 64--80.

\bibitem[{Morato et~al.(2023)Morato, Normey-Rico, and Sename}]{Morato23}
Morato, M.M., Normey-Rico, J.E., and Sename, O. (2023).
\newblock Sufficient conditions for convergent recursive extrapolation of qlpv scheduling parameters along a prediction horizon.
\newblock \emph{IEEE Transactions on Automatic Control}, 68(6), 3182--3193.

\bibitem[{Nezami et~al.(2023)Nezami, Karachalios, Schildbach, and Abbas}]{Maryam}
Nezami, M., Karachalios, D.S., Schildbach, G., and Abbas, H.S. (2023).
\newblock On the design of nonlinear {MPC} and {LPVMPC} for obstacle avoidance in autonomous driving*.
\newblock In \emph{2023 9th International Conference on Control, Decision and Information Technologies (CoDIT)}, 1--6.

\bibitem[{Nguyen-Van and Hori(2014)}]{VaHo14}
Nguyen-Van, T. and Hori, N. (2014).
\newblock Discretization of nonautonomous nonlinear systems based on continualization of an exact discrete-time model.
\newblock \emph{Journal of Dynamic Systems, Measurement, and Control}, 136(2), 1 -- 9.

\bibitem[{Nocedal and Wright(2006)}]{NoceWrig06}
Nocedal, J. and Wright, S.J. (2006).
\newblock \emph{Numerical Optimization}.
\newblock Springer, New York, NY, USA, 2e edition.

\bibitem[{Qin and Badgwell(2003)}]{QiBa03}
Qin, S. and Badgwell, T. (2003).
\newblock A survey of industrial model predictive control technology.
\newblock \emph{Control Engineering Practice}, 11, 733--746.

\bibitem[{Rawlings and Mayne(2017)}]{RaMaDi17chapter8}
Rawlings, J.B. and Mayne, D.Q. (2017).
\newblock \emph{Model Predictive Control: Theory, Computation, and Design}, chapter~8.
\newblock Nob Hill Publishing, LLC.
\newblock Numerical Optimal Control.

\bibitem[{T{\'o}th(2010)}]{To10}
T{\'o}th, R. (2010).
\newblock \emph{Modeling and Identification of Linear Parameter-Varying Systems}.
\newblock Springer--Verlag, Berlin, Heidelberg.

\bibitem[{W\"{a}chter and Biegler(2006)}]{WaBi06}
W\"{a}chter, A. and Biegler, L. (2006).
\newblock On the implementation of an interior-point filter line-search algorithm for large-scale nonlinear programming.
\newblock \emph{Math. Program.}, 106, 25–57.

\end{thebibliography}
                                                   







\appendix
\section{} \label{a:NonCondMat}   

According to the constraints on the state and inputs, the linear inequality constraints in \eqref{eq:opt1}
are defined by
\begin{align}
    \begin{aligned}
    \label{eq:constraints}
      &\mathcal{X} = \{ x \in \mathbb{R}^{n_{{\mathrm{x}}}}~|~G_\text{x} x \leq h_\text{x} \},\\
       &\mathcal{U} = \{ u \in \mathbb{R}^{n_u}~|~G_\text{x} u \leq h_\text{u} \},\\
    \end{aligned} 
\end{align}
where $G_\text{x}, G_\text{u}$ 
are matrices and  $h_\text{x}, h_\text{u}$ 
are vectors, all with appropriate dimensions,   such that the inequalities represent polyhedral sets.

The matrices $C(\bar{p}), G^\text{z}, h^\text{z}$ in \eqref{e:lpvmpcsim}  are given, respectively, as follows
\begin{align*}
    C(\bar{p})&=\left[\begin{array}{ccc|cccc}
        B(p_{0}) & 0  & \cdots  & -\mathbb{I}_{n} & 0 & 0 &  \cdots \\
        0 & B(p_{1})  & \cdots  &  A(p_{1}) & -\mathbb{I}_{n}  & 0 & \cdots\\
        \vdots  & \ddots & \ddots  & \vdots & \ddots & \ddots & \ddots  \\
    \end{array}\right]
\\
    {G}^{\rm z}&=\text{blkdiag}(\bar{G}_{\rm u},\bar{G}_{\rm x}),\quad
    {h}^{\rm z}=\bbma \bar{h}_{\rm u} \\ \bar{h}_{\rm x}\ebma,
\end{align*}
where
\begin{align*}
\bar{G}_{\rm u}&=\mathbb{I}_{N}\otimes G_\text{u},&& \bar{G}_{\rm x}=\text{blkdiag}(\mathbb{I}_{N-1}\otimes G_\text{x},G_\text{xf}),\\
 \bar{h}_{\rm u}&=1_{N}\otimes h_\text{u},&& \bar{h}_{\rm x}=\bbma
 1_{{N-1}}\otimes h_\text{x} \\
          h_\text{xf} \ebma
\end{align*}
and  $1_{N}\in\R^N$, $1_{N-1}\R^{N-1}$ denote vectors with all entries equal to one.

The matrices $\bar{A}(\bar{p})$ and $\bar{B}(\bar{p})$ in \eqref{e:AugLPVSys} are given by the following expressions
\begin{align*}
     \bar{A}(\bar{p})\!&=\!\left[\begin{matrix}
     A(p_{0})\\A(p_{1})A(p_{0})\\\vdots\\\prod\limits_{k=0}^{N-1} A(p_i)
 \end{matrix}\right],\\
 \bar{B}(\bar{p})\!&=\!\!\left[\begin{matrix} B(p_{0}) & 0 & \cdots & 0\\ A(p_{1})B(p_{0}) & B(p_{1}) & \cdots & 0\\ \vdots & \vdots & \ddots & \vdots\\\prod\limits_{k=1}^{N-1}\!\!\! A(p_{i})B(p_0) & \prod\limits_{k=2}^{N-1}\!\!\!A(p_{i})B(p_1) & \cdots & B(p_{N-1}) \end{matrix}\right]
\end{align*}
Finally, ${G}^{\rm u}(\bar{p})$, ${h}^{\rm u}(\bar{p})$ in (\ref{e:LPVMPCQP}b) are given as follows:
\[
     {G}^{\rm u}(\bar{p})=\bbma \bar{G}_{\rm u} \\\bar{G}_{\rm x}\bar{B}(\bar{p}) \ebma,\quad
       {h}^{\rm u}(\bar{p})=\bbma \bar{h}_{\rm u}\\ \bar{h}_{\rm x}-\bar{G}_{\rm x}\bar{A}(\bar{p})x_0\ebma. 
\]

\end{document}